\input amstex.tex

\documentstyle{amsppt}
\magnification \magstephalf
\pagewidth{5.4in} \pageheight{7.7in} \vcorrection{-0.2in}
\NoRunningHeads
\def\makefootline{\baselineskip=49pt \line{\the\footline}}

\define\ep{\varepsilon}

\define\th{\theta}
\redefine\Im{\operatorname{Im}}

\define\trace{\operatorname{trace}}
\NoBlackBoxes \magnification\magstep1

\topmatter
\title
A negative mass theorem for surfaces of positive genus
\endtitle
\author  K. Okikiolu
\endauthor
\thanks{Revised October 11, 2008.   MSC classes: 58J50,  35J60,
53C20.\newline The author would like to acknowledge the support of
the Institute for Advanced Study.} \endthanks \abstract Let $M$ be a
closed surface. For a metric $g$ on $M$,
 denote the Laplace-Beltrami operator by
$\Delta=\Delta_g$.  We define $\trace\Delta^{-1}=\int_M m(p)\,dA$,
where $dA$ is the area element for $g$ and $m(p)$ is the Robin
constant  at the point $p\in M$, that is the value of the Green
function $G(p,q)$ at $q=p$ after the logarithmic singularity has
been subtracted off. Since $\trace\Delta^{-1}$ can also be obtained
by regularization of the spectral zeta function, it is  a spectral
invariant. Heuristically it represents the sum of squares of the
wavelengths of the surface.  We define the $\Delta$-mass of $(M,g)$
to equal $(\trace\Delta_g^{-1}-\trace \Delta_{S^2,A}^{-1})/A$, where
$\Delta_{S^2,A}$ is the Laplacian on the round sphere of area $A$.
This is an analog for closed surfaces of the ADM mass from general
relativity.
 We show that if $M$ has positive genus, the minimum
of the $\Delta$-mass on each conformal class
is negative and attained by a smooth metric.
For this minimizing metric, there is a sharp logarithmic Hardy-Littlewood-Sobolev inequality and a Moser-Trudinger-Onofri type inequality.
\newline
\noindent
\endabstract
\endtopmatter

\noindent{\bf Section 1. Introduction}.
\medskip

Let $M$ be a closed Riemann surface and let $g$ be a metric on $M$
compatible with the complex structure.  With respect to complex
coordinates $z$, the metric $g$ is the real part of the K\"{a}hler
metric
$$
e^u \, dz\otimes d\bar z,
$$
for some smooth real valued function $u$ on $M$.  The area element
is
$$
dA\ =\ \frac{ie^u }2 dz\wedge d\bar z.
$$
Denote the total area by $A$. The Laplace Beltrami operator for the
metric $g$ is given by
$$
\Delta\ =\ \Delta_g\ =\ -4e^{-u}\partial_z\partial_{\bar z}.
$$
(This is sometimes called {\it geometer's Laplacian}: note the
sign.) The Green's function for the metric $g$ is the smooth real
valued function $G$ on $M\times M\setminus\{(p,p):p\in M\}$ such
that
$$
\int_M G(p,q)\ \Delta f(q)\,dA(q)\ =\ f(p)\ -\ \frac1A \int_M f\,dA
$$
for smooth functions $f$ on $M$. It follows that $G$ is symmetric,
$\int G(p,q)\,dA(q)=0$, and
$$
\Delta_p G(p,q)\ =\ -\frac1A\qquad\qquad\text{when }p\neq q. \tag 1.1
$$
For the smooth function $f$ on $M$ we define
$$
\Delta^{-1}f(p)\ =\ \int_M G(p,q)f(q)\,dA(q).
$$
If $d(p,q)$ is the geodesic distance between $p$ and $q$ in the metric $g$
then there exists a smooth function $m$ on $M$ such that
$$
G(p,q)\ =\ \frac1{2\pi}\log d(p,q)\ +\ m(p)\ +\ O(d(p,q)),\qquad
\text{ as }d(p,q)\to 0.
$$
The value $m(p)$ is known as the {\it Robin constant} at $p$. We
define
$$
\trace \Delta^{-1}_g \ =\ \int_M m_g\,dA.
$$
 This is  a spectral invariant for $\Delta$, since
 it can be obtained from the spectral zeta function associated to $\Delta$, see [S1],
 [S2], [M3], or [Ok1]. Heuristically it represents
the {\it sum of squares of the wavelengths} of the surface (up to a
 constant). It is convenient to
normalize to get a scale invariant quantity.  Indeed, define the
$\Delta$-mass to be
$$
\Cal M(g)\ =\ \frac{\trace\Delta_g^{-1}-\trace\Delta_{S^2,A}^{-1}}{A},
$$
where $\Delta_{S^2,A}$ is the Laplacian for the round metric on
$S^2$ with area $A$.  Then the above results show that $\Cal M(g)$
is always positive when $g$ is a metric on $S^2$.   In this paper we
show the following:

\proclaim{Theorem 1} (Negative mass theorem for positive genus
surfaces.) Given a metric $g$ on the closed surface $M$ of positive
genus,   there exists a conformal metric $e^\phi g$ such that $\Cal
M( e^\phi g)<0$. In fact $e^\phi g$ can be chosen to minimize $\Cal
M$ within the conformal class.
\endproclaim

When $M$ is a torus, this was proved in [Ok2].  We remark that
Theorem 1 fails on the sphere.  This follows from the logarithmic
Hardy-Littlewood-Sobolev inequality for the sphere [On], [CL], [B]:

\proclaim{Theorem (Morpurgo [M2])}  For a metric $g$ on the sphere
$S^2$, the value $\Cal M(g)$ is strictly positive unless $g$ is
round.
\endproclaim

From [Ok1], we immediately obtain the following corollary to Theorem
1.\medskip

\proclaim{Corollary 2} (Analogs of Logarithmic HLS inequality and
the Moser-Trudinger-Onofri Inequality for general surfaces.) If the
metric $g$ minimizes $\Cal M$ within its conformal class, then
$$
 \frac1{4\pi}\int_{M} \psi\,e^\psi \,dA\ -\ \frac1A\int_{M}\, e^\psi \Delta^{-1}
 e^\psi \,dA \ \geq\ 0
$$
for all  functions $\psi:M\to\Bbb R$ with $ \int_{M} e^\psi \,dA=A$
such that $\int_{M} \psi\,e^\psi \,dA$ is finite.  Here, $dA$ and
$\Delta$ are associated to $g$. Moreover, for $\psi\in C^\infty(M)$,
$$
\frac{1}{16\pi}\int_M \psi \Delta \psi\,dA\ -\ \log\left(
\frac1A\int_M e^{\psi}\,dA\right)\ +\ \frac1A\int_M \psi\,dA\ \geq\
0.
$$
\endproclaim
For some related results, see [Ch], [CheC], [DJLW], [LL1], [LL2],
[M2], [M3], [NT], [Ok1], [Ok2], [OPS], [S2].

\medskip

 \remark{Remark}  The quantity $\Cal M(g)$  can be viewed as an analog of the ADM mass.
 Indeed, writing $K(p)$ for the Gaussian curvature of
$g$ at $p$, it is shown in [S1], [S2],  that for any metric $g$ on
the $2$-sphere, the natural analog of the ADM mass for metrics on
the sphere is the constant
$$
m_g(p)\ -\ \frac1{2\pi}\Delta^{-1}K(p)\ =\ \frac1A
\trace\Delta^{-1}_g. \tag 1.2
$$
Although the left hand side of (1.2) is not constant in general for
surfaces of higher genus, the right hand  can be thought of as the
analog of the mass.
 (The left hand side of (1.2) is constant  for the {\it canonical
metric}, a fact  we use in next section.)  For a probabilistic
interpretation of $\trace\Delta^{-1}$, see [DS1].  There it is shown
that $\trace\Delta^{-1}$ is the constant term in an asymptotic
expansion in $\ep$, of the time it takes a Brownian particle
starting at a randomly chosen point on the surface to get
$\ep$-close to another randomly chosen point.
\endremark
\medskip

\noindent{\bf Section 2. The Proof.}
\medskip

There are two main ingredients in the proof of this result.  The
first is an identity concerning the Arakelov Green's function which
is used in the construction of the Arakelov metric.  The second is a
delicate result on the mean field equation on surfaces proved in
[DJLW]. That paper gives  conditions under which a general mean
field equation has a solution.  Here we show that for the particular
case of the {\it canonical metric} on $M$ and the mean field
equation arising from $\trace\Delta^{-1}$,  the conditions of the
[DJLW] theorem are satisfied. We start by recalling the way that
Robin's constant and the sum of squares of the wavelengths change
under a conformal change of the metric.

\proclaim{Proposition 2.1: Conformal change of the Robin constant}
If $\phi$ is a smooth function on $M$ then
$$
m_{e^\phi g}(p)\ =\  m _g(p)\ +\ \frac{\phi}{4\pi}\ -\  \frac{2}{ A_\phi}
(\Delta^{-1}_g e^\phi)(p) \ +\  \frac1{A_\phi^2}\int_M e^\phi
\Delta_g^{-1}e^\phi\,dA,
$$
where
$$
A_\phi\ =\ \int_M e^{\phi}\,dA.
$$
\endproclaim
\noindent For the proof, see for example [S1], [S2], [M3] or [Ok1].

\proclaim{Proposition 2.2: Conformal change of $\trace\Delta^{-1}$
(Morpurgo's Formula)}  If $\phi$ is a smooth function on $M$, then
$$
\trace \Delta^{-1}_{e^\phi g}\ =\ \int_M m_g e^\phi \,dA\ +\ \frac1{4\pi}\int_M
\phi\, e^\phi \,dA\ -\ \frac1{A_{\phi}} \int_M e^\phi \Delta^{-1}_g e^\phi\,dA.
\tag 2.1
$$
\endproclaim

Now we discuss the Canonical metric and the Arakelov Green's
function. We refer the reader to [W] and [F] for more details on
this subject. If $M$ is a Riemann surface of genus $H$, there exists
a metric $g$ on $M$ known as the {\it canonical metric} which is
compatible with the complex structure. It is defined by taking the
Jacobian embedding of the Riemann surface $M$ into
 a $2H$-dimensional torus, and pulling back the flat metric
on the torus to $M$.  Indeed,   Let $\{A_j,B_j\}$ be a symplectic
homology basis for $H_1(M,\Bbb Z)$ satisfying the intersection
pairings
$$
\#[A_i,A_j]=0,\qquad\#[B_i,B_j]=0,\qquad \#[A_i,B_j]=\delta_{ij}.
$$
Take a basis $\th_j$ for the space of homomorphic $1$-forms
satisfying
$$
\int_{A_j}\th_k\ =\ \delta_{jk}.
$$
Then the  period matrix $\Omega_{ij}$  given by
$$
\Omega_{ij}\ =\ \int_{B_i}\th_j
$$
is positive definite. The Jacobian variety associated to $M$ is
$$
J(M)\ =\ \Bbb C^H/(\Bbb Z^h+\Omega\Bbb Z^H).
$$
The {\it Abel map} gives an embedding of $M$ into $J(M)$,
$$
I:z\to \int_{z_0}^z (\th_1,\dots,\th_H).
$$
 The canonical K\"{a}hler metric on
$M$ is given in terms of local holomorphic coordinates $z$ by
$$
\mu(z)\ dz\otimes d\bar z,\qquad\text{where}\qquad
\mu(z)\ =\ \frac1{H}\left(\sum_{j,k=1}^H (\Im \Omega)^{-1}_{jk}\frac{d\th_j}{dz}
\frac{d\bar\th_k}{d\bar z}\right).
$$
The real part of the K\"{a}hler metric is the Riemannian metric $g$.
It can be checked that this metric has unit area. The Green's
function for this metric is known as the {\it Arakelov Green's
function} and the following result is well known.

\proclaim{Proposition 2.3}  If $M$ is a closed Riemann surface of
genus $H$ and if $g$ is the canonical metric on $M$ with unit area,
then the Robin constant $m(p)$ for $g$ satisfies
$$
\Delta m(p)\ =\ 2H\ -\ 2\ + \frac{K(p)}{2\pi}.
$$
\endproclaim

\noindent{\bf Proof of Proposition 2.3}.
 The Gaussian curvature $K$ is given by
$$
2\partial_z\partial_{\bar z}\log \mu\ =\ -\mu K.
$$
The Arakelov Green's function is given by
$$
\multline
G(z,w)\ =\ -\frac1{2\pi}\log|E(z,w)|\ +\ \frac12\sum_{j,k=1}^H  (\Im \Omega)^{-1}_{jk}\Im(Z-W)_j \Im(Z-W)_k
\\ +\ \frac{m(z)}2-\ \frac{\log \mu(z)}{8\pi}\ +\ \frac{m(w)}2\
\ -\ \frac{\log \mu(w)}{8\pi}.
\endmultline
$$
Here,  $Z=I(z)$ and $W=I(w)$  and $E(z,w)$ is the {\it prime form}
which plays the role of $z-w$, is holomorphic  in $z$ and $w$, and
transforms as a $(-1/2,-1/2)$ form in each variable.  We notice that
$\log |E(z,w)|$ is harmonic, and so from (1.1) we have that for
$w\neq z$,
$$
\mu(z) \ =\ 4\partial_z\partial_{\bar z}G(z,w)\ =\ H\mu(z)\ +\ 2\partial_z\partial_{\bar z}m(z)
\ +\ \frac{\mu(z) K(z)}{4\pi}.
$$
Hence we see that
$$
4\mu^{-1}\partial_z\partial_{\bar z}m\ =\ 2\ -\ 2H\ -\ \frac{K}{2\pi}.\qquad\qquad\qed
$$
\medskip\medskip

Theorem 1 is now an application of the following result on the mean
field equation which is obtained from  Theorem 1.2 of [DJLW] and its
proof.
\medskip

\proclaim{Theorem 2.4. [DJLW]} Let $(M,g)$ be a closed surface of
unit area and let $h$ be a smooth positive function on $M$.  Suppose
$p_0$ is a point at which $8\pi m+2\log h$ attains its maximum
value, and suppose in addition that
$$
\Delta \log h(p_0)\ <\ 8\pi \ -\ 2K(p_0).
$$
Then the minimum of the functional
$$
J(u)\ =\ \frac1{16\pi}\int_M|\nabla u|^2\,dA\ +\  \int_M u\,dA\ -\
 \log\int_M he^u\,dA \tag 2.2
$$
over functions $u$ in the Sobolev space $H^1(M)$ is attained at a
smooth function $u$ satisfying
$$
\Delta u\ =\ 8\pi he^u \ -\ 8\pi. \tag 2.3
$$
Moreover, for this minimum point $u$ we have
$$
 J(u)\ <\ -\left( 1 +\log \pi \ +\ \max_{p\in M}
(4\pi m_g(p)+\log h(p))\right).\tag 2.4
$$
\endproclaim

\noindent{\bf Proof of Theorem 1}.  We take $g$ to be the canonical
metric on $M$, and we set
$$
 h=e^{-4\pi m_g}.
$$
Then by Proposition 2.3, we have
$$
\Delta \log h\ =\ -4\pi \Delta m_g\ =\ 8\pi\ -\ 8\pi H\ -\ 2K\ <\ 8\pi-2K.
$$
Hence we obtain the conclusion of Theorem 2.4.  From (2.3), the
function $u$ satisfies
$$
\int_M he^u\,dA\ =\ 1.
$$
from (2.2) and (2.3) we see that
$$
J(u)\ =\ \frac12\int_M u(he^u+1)\,dA.
$$
However, writing
$$
\phi\ =\ u-4\pi m_g,
$$
we have
$$
\int_M e^\phi\,dA\ =\ 1,
$$
and
$$
\Delta_g^{-1} e^\phi\ =\ \frac1{8\pi}\left( u-\int_M u\right). \tag 2.5
$$
Hence from (2.1),
$$
\trace \Delta^{-1}_{e^\phi g}\ =\ \frac1{8\pi} \int_M u(he^u+1)
\ =\ \frac{J(u)}{4\pi}.
$$
Now from (2.4) and the fact that
$$
\trace\Delta^{-1}_{S^2,1}\ =\ \frac{-1-\log \pi}{4\pi},
$$
we see that
$$
\trace\Delta^{-1}_{e^\phi g}\ <\ \trace\Delta^{-1}_{S^2,1}.
$$
In fact we remark that $e^\phi g$ minimizes $\trace\Delta^{-1}$
among unit area  metrics  in the conformal class of $g$. Indeed,
from [Ok] Theorem 1,  we  conclude  that the minimum of
$\trace\Delta_{e^\psi g}^{-1}$ among conformal factors with $\int
e^\psi=1$, must in fact be attained at a metric $e^\psi g$ with
$\int_M e^\psi\,dA=1$, which must also satisfy the Euler-Lagrange
equation, namely that the Robin constant $m_{e^\psi g}(p)$ is
constant:
$$
\Delta_g^{-1} e^\psi\ =\ \frac1{8\pi}\left( \psi\ +\ 4\pi m_g\ -\ \int_M (\psi+4\pi m_g)\,dA\right).
$$
However, setting $v=\psi+4\pi m_g$, we find that $v$ is a critical
metric for $J$, and hence
$$
\trace\Delta^{-1}_{e^\psi g}\ =\ \frac{J(v)}{4\pi}\ \geq\  \frac{J(u)}{4\pi}\ =\ \trace\Delta^{-1}_{e^\phi g}.
$$
We also remark that  in general the metric $e^\phi g$ need not
coincide with  the canonical metric or the constant curvature metric
or the Arakelov metric, see [Ok2]. On a long thin rectangular torus,
the minimizer  is close to being a round sphere  with a short worm
hole joining the poles.

\medskip
I would like to thank Richard Wentworth for helpful discussions.

\medskip
\centerline{References}
 \roster

\item"[B]" Beckner, W.: Sharp Sobolev inequalities on the sphere and the
Moser-Trudinger inequality. {\it Annals of Math.} {\bf 138} (1993),
213-242.

\item"[CL]" Carlen, E.,  Loss, M.: Competing symmetries, the logarithmic HLS inequality
and Onofri's inequality on $S^n$.  {\it Geometric and Functional
Analysis} {\bf 2} (1992) 90--104.

\item"[Ch]" Chang, S.-Y. A.: Conformal invariants and partial differential equations.
{\it Bull. Amer. Math. Soc.} {\bf  42} (2005),  365--393.

\item"[CheC]" Chen, C.-C., Lin, C.-S.: Sharp estimates for solutions
of multi-bubbles in compact Riemann surfaces. Comm. Pure Appl. Math.
{\bf 55}  (2002),  no. 6, 728--771.

\item"[DJLW]" Ding, W., Jost, J., Li, J., Wang, G.: The differential
equation $\Delta u=8\pi -8\pi he^u$ on a compact Riemann surface.
Asian J. Math. {\bf 1}, 230--248 (1997).

\item"[DS1]"  Doyle, P., Steiner, J.: Spectral invariants and playing hide and
seek on surfaces.
{\it Preprint}.

\item"[DS2]"  Doyle, P.,  Steiner, J.: Blowing bubbles on the torus.
{\it Preprint}.

\item"[F]" Fay, T.:  Theta functions on Riemann surfaces. Ann. Math.
{\bf 119}, 387 (1994).

\item"[LL1]" Lin, C.-S., Lucia, M.: Uniqueness of solutions for a mean
field equation on the torus. J. Diff. Eqs. {\bf 229}(1), 172--185 (2006).

\item"[LL2]" Lin, C.-S., Lucia, M.: One-dimensional symmetry of periodic minimizers for a mean field equation.
Ann. Sc. Norm. Super. Pisa Cl. Sci. (5) {\bf 6}(2), 269–290 (2007)

\item"[M1]"  Morpurgo, C.:  The logarithmic Hardy-Littlewood-Sobolev inequality
and extremals of zeta functions on $S^n$. {\it Geom. Funct. Anal.}
{\bf 6} (1996),  146--171.

\item"[M2]"  Morpurgo, C.: Zeta functions on $S\sp 2$. Extremal Riemann
surfaces (San Francisco, 1995), 213--225, {\it Contemp. Math.}, {\bf
201},

\item"[M3]" Morpurgo, C.:
 Sharp inequalities for functional integrals and traces of conformally invariant
 operators.
 {\it Duke Math. J.} {\bf 114} (2002),  477--553.

\item"[NT]" Nolasco, M., Tarantello, G.: On a sharp Sobolev-type inequality on two-dimensional compact
manifolds. Arch. Ration. Mech. Anal. {\bf 145}, 161–195 (1998)

\item"[Ok1]" K. Okikiolu:  Extremals for Logarithmic HLS inequalities on compact
manifolds.  GAFA {\bf 107}(5),
1655–1684 (2008)

\item"[Ok2]" Okikiolu, K: A negative mass theorem for the 2-torus.
{\it Comm. Math. Physics}, (To Appear.)

\item"[On]" Onofri, E.: On the positivity of the effective action in a theory of
random surfaces.
{\it Comm. Math. Phys.} {\bf  86} (1982),  321--326.

\item"[OPS]"  Osgood, B.,  Phillips, R.,   Sarnak, P.: Extremals of determinants of
 Laplacians.
{\it J. Funct. Anal.} {\bf 80} (1988),  148--211.

\item"[S1]"  Steiner, J: {\it Green's Functions, Spectral Invariants, and a Positive Mass on Spheres}.
Ph. D. Dissertation, University of California San Diego, June 2003.

\item"[S2]" Steiner, J:  A geometrical mass and its extremal properties for metrics on $S\sp 2$.
{\it Duke Math. J.}  {\bf 129} (2005), 63--86.

\item"[W]" Wentworth, R.:  The asymptotics of the Arakelov-Green's
function and Faltings' delta invariant.  {\it Comm. Math. Phys.}
{\bf 137}, 427--459.  (1991)

\endroster
\medskip\medskip

\centerline{Kate Okikiolu}

\centerline{University of California, San Diego}

\centerline{okikiolu\@math.ucsd.edu}

\end